\input amstex
\documentstyle{amsppt}
%
\catcode`@=11
\redefine\output@{%
  \def\break{\penalty-\@M}\let\par\endgraf
  \ifodd\pageno\global\hoffset=105pt\else\global\hoffset=8pt\fi  
  \shipout\vbox{%
    \ifplain@
      \let\makeheadline\relax \let\makefootline\relax
    \else
      \iffirstpage@ \global\firstpage@false
        \let\rightheadline\frheadline
        \let\leftheadline\flheadline
      \else
        \ifrunheads@ 
        \else \let\makeheadline\relax
        \fi
      \fi
    \fi
    \makeheadline \pagebody \makefootline}%
  \advancepageno \ifnum\outputpenalty>-\@MM\else\dosupereject\fi
}
\def\Beta{\mathchar"0\hexnumber@\rmfam 42}
\catcode`\@=\active
\nopagenumbers
\chardef\textvolna='176
\def\negskp{\hskip -2pt}
\def\rank{\operatorname{rank}}
\def\Sym{\operatorname{Sym}}
\def\blue#1{#1}

\catcode`#=11\def\diez{#}\catcode`#=6
\catcode`&=11\catcode`&=4
\catcode`_=11\def\podcherkivanie{_}\catcode`_=8
\catcode`~=11\def\volna{~}\catcode`~=\active
\def\mycite#1{\cite{\blue{#1}}\immediate\special{ps:
     ShrHPSdict begin /ShrBORDERthickness 0 def}}
\def\myciterange#1#2#3#4{\cite{\blue{#2#3#4}}\immediate\special{ps:
     ShrHPSdict begin /ShrBORDERthickness 0 def}}
\def\mytag#1{%
    \tag#1}
\def\mythetag#1{\thetag{\blue{#1}}\immediate\special{ps:
     ShrHPSdict begin /ShrBORDERthickness 0 def}}
\def\myrefno#1{\no#1}
\def\myhref#1#2{\blue{#2}\immediate\special{ps:
     ShrHPSdict begin /ShrBORDERthickness 0 def}}
\def\myEarXivlink{\myhref{http://arXiv.org}{http:/\negskp/arXiv.org}}

\def\mytheorem#1{\csname proclaim\endcsname{Theorem #1}}
\def\mytheoremwithtitle#1#2{\csname proclaim\endcsname{Theorem #1#2}}
\def\mythetheorem#1{\blue{#1}\immediate\special{ps:
     ShrHPSdict begin /ShrBORDERthickness 0 def}}
\def\mylemma#1{\csname proclaim\endcsname{Lemma #1}}
\def\mylemmawithtitle#1#2{\csname proclaim\endcsname{Lemma #1#2}}

\def\mycorollary#1{\csname proclaim\endcsname{Corollary #1}}

\def\mydefinition#1{\definition{Definition #1}}
\def\mythedefinition#1{\blue{#1}\immediate\special{ps:
     ShrHPSdict begin /ShrBORDERthickness 0 def}}
\def\myconjecture#1{\csname proclaim\endcsname{Conjecture #1}}
\def\myconjecturewithtitle#1#2{\csname proclaim\endcsname{Conjecture #1#2}}


\pagewidth{360pt}
\pageheight{606pt}
\topmatter
\title
On the equivalence of cuboid equations\\ and their factor equations. 
\endtitle
\rightheadtext{On the equivalence of cuboid equations \dots}
\author
Ruslan Sharipov
\endauthor
\address Bashkir State University, 32 Zaki Validi street, 450074 Ufa, Russia
\endaddress
\email\myhref{mailto:r-sharipov\@mail.ru}{r-sharipov\@mail.ru}
\endemail
\abstract
    An Euler cuboid is a rectangular parallelepiped with integer edges and 
integer face diagonals. An Euler cuboid is called perfect if its space diagonal 
is also integer. Some Euler cuboids are already discovered. As for perfect cuboids, 
none of them is currently known and their non-existence is not yet proved. Euler 
cuboids and perfect cuboids are described by certain systems of Diophantine equations.
These equations possess an intrinsic $S_3$ symmetry. Recently they were factorized
with respect to this $S_3$ symmetry and the factor equations were derived. In the
present paper the factor equations are shown to be equivalent to the original 
cuboid equations regarding the search for perfect cuboids and in selecting Euler
cuboids. 
\endabstract
\subjclassyear{2000}
\subjclass 11D41, 11D72, 13A50, 13F20\endsubjclass
\endtopmatter
\TagsOnRight
\document

\head
1. Introduction.
\endhead
     The search for perfect cuboids has the long history. The reader can follow this
history since 1719 in the references \myciterange{1}{1}{--}{44}. In order to write
the cuboidal Diophantine equations we use the following polynomials:
$$
\xalignat 2
&\hskip -2em
p_{\kern 1pt 0}=x_1^2+x_2^2+x_3^2-L^2,
&&p_{\kern 1pt 1}=x_2^2+x_3^2-d_1^{\kern 1pt 2},\\
\vspace{-1.7ex}
\mytag{1.1}\\
\vspace{-1.7ex}
&\hskip -2em
p_{\kern 1pt 2}=x_3^2+x_1^2-d_2^{\kern 1pt 2},
&&p_{\kern 1pt 3}=x_1^2+x_2^2-d_3^{\kern 1pt 2}.
\endxalignat
$$
Here $x_1$, $x_2$, $x_3$ are edges of a cuboid, $d_1$, $d_2$, $d_3$ are its face 
diagonals, and $L$ is its space diagonal. An Euler cuboid is described by a system
of three Diophantine equations. In terms of the polynomials \mythetag{1.1} these 
equations are written as
$$
\xalignat 3
&\hskip -2em
p_{\kern 1pt 1}=0,&&p_{\kern 1pt 2}=0,&&p_{\kern 1pt 3}=0.
\mytag{1.2}
\endxalignat
$$ 
In the case of a perfect cuboid the number of equations is greater by one, i\.\,e\.
instead of the equations \mythetag{1.2} we write the following system of four
equations:
$$
\xalignat 4
&\hskip -2em
p_{\kern 1pt 0}=0,&&p_{\kern 1pt 1}=0,&&p_{\kern 1pt 2}=0,&&p_{\kern 1pt 3}=0.
\mytag{1.3}
\endxalignat
$$
\par
     The permutation group $S_3$ acts upon the cuboid variables $x_1$, $x_2$, $x_3$,
$d_1$, $d_2$, $d_3$, and $L$ according to the rules expressed by the formulas
$$
\xalignat 3
&\hskip -2em
\sigma(x_i)=x_{\sigma i},
&&\sigma(d_i)=d_{\sigma i},
&&\sigma(L)=L.
\mytag{1.4}
\endxalignat
$$ 
The variables $x_1$, $x_2$, $x_3$ and $d_1$, $d_2$, $d_3$ are usually arranged into
a matrix:
$$
\hskip -2em
M=\Vmatrix x_1 & x_2 &x_3\\
\vspace{1ex}
d_1 & d_2 & d_3\endVmatrix.
\mytag{1.5}
$$
The rules \mythetag{1.4} means that $S_3$ acts upon the matrix \mythetag{1.5} by
permuting its columns. Applying the rules \mythetag{1.4} to the polynomials
\mythetag{1.1}, we derive 
$$
\xalignat 2
&\hskip -2em
\sigma(p_i)=p_{\kern 0.5pt \sigma i},
&&\sigma(p_{\kern 1pt 0})=p_{\kern 1pt 0}.
\mytag{1.6}
\endxalignat
$$
\par
     The polynomials $p_{\kern 1pt 0}$, $p_{\kern 1pt 1}$, $p_{\kern 1pt 2}$, 
$p_{\kern 1pt 3}$ in \mythetag{1.1} are treated as elements of the polynomial
ring $\Bbb Q[x_1,x_2,x_3,d_1,d_2,d_3,L]$. For the sake of brevity we denote 
$$
\hskip -2em
\Bbb Q[x_1,x_2,x_3,d_1,d_2,d_3,L]=\Bbb Q[M,L],
\mytag{1.7}
$$
where $M$ is the matrix given by the formula \mythetag{1.5}. 
\mydefinition{1.1} A polynomial $p\in\Bbb Q[M,L]$ is called multisymmetric if it
is invariant with respect to the action \mythetag{1.4} of the group $S_3$. 
\enddefinition
     Multisymmetric polynomials constitute a subring in the ring \mythetag{1.7}.
We denote this subring through $\Sym\!\Bbb Q[M,L]$. The formulas \mythetag{1.6}
show that the polynomial $p_{\kern 1pt 0}$ belongs to the subring $\Sym\!\Bbb 
Q[M,L]$, i\.\,e\. it is multisymmetric, while the polynomials $p_{\kern 1pt 1}$, 
$p_{\kern 1pt 2}$, $p_{\kern 1pt 3}$ are not multisymmetric. Nevertheless, the 
system of equations \mythetag{1.2} in whole is invariant with respect to the 
action of the group $S_3$. The same is true for the system of equations 
\mythetag{1.3}.\par
     The polynomials $p_{\kern 1pt 1}$, $p_{\kern 1pt 2}$, $p_{\kern 1pt 3}$
generate an ideal in the ring $\Bbb Q[M,L]$. It is natural to call it the 
{\it cuboid ideal\/} and denote this ideal through
$$
\hskip -2em
I_{\text{C}}=\bigl<p_{\kern 1pt 1},p_{\kern 1pt 2},p_{\kern 1pt 3}\bigr>.
\mytag{1.8}
$$
Similarly, one can define the {\it perfect cuboid ideal} 
$$
\hskip -2em
I_{\text{PC}}=\bigl<p_{\kern 1pt 0},p_{\kern 1pt 1},p_{\kern 1pt 2},
p_{\kern 1pt 3}\bigr>.
\mytag{1.9}
$$
Each polynomial equation $p=0$ with $p\in I_{\text{C}}$ follows from the 
equations \mythetag{1.2}. Therefore such an equation is called a {\it cuboid
equation}. Similarly, each polynomial equation $p=0$ with $p\in I_{\text{PC}}$
follows from the equations \mythetag{1.3}. Such an equation is called a {\it
perfect cuboid equation}.\par
     The symmetry approach to the equations \mythetag{1.2} and \mythetag{1.3}
initiated in \mycite{45} leads to studying the following ideals in the ring
of multisymmetric polynomials:
$$
\xalignat 2
&\hskip -2em
I_{\text{C\kern -0.7pt\_\kern 0.5pt sym}}=I_{\text{C}}\cap\Sym\!\Bbb Q[M,L],
&&I_{\text{PC\kern -0.7pt\_\kern 0.5pt sym}}=I_{\text{PC}}\cap\Sym\!\Bbb Q[M,L].
\quad
\mytag{1.10}
\endxalignat
$$
\mydefinition{1.2} A polynomial equation of the form $p=0$ with $p\in I_{\text{C
\kern -0.7pt\_\kern 0.5pt sym}}$ or with $p\in I_{\text{PC\kern -0.7pt\_\kern 0.5pt 
sym}}$ is called an {\it $S_3$ factor equation} for the Euler cuboid equations 
\mythetag{1.2} or for the perfect cuboid equations \mythetag{1.3} respectively. 
\enddefinition
     The ideal $I_{\text{PC\kern -0.7pt\_\kern 0.5pt sym}}$ from \mythetag{1.10}
was studied in \mycite{46} (there it was denoted through $I_{\text{sym}}$). The 
polynomial $p_{\kern 1pt 0}$ used as a generator in \mythetag{1.9} is multisymmetric 
in the sense of the definition~\mythedefinition{1.1}. Therefore it belongs to the 
ideal $I_{\text{PC\kern -0.7pt\_\kern 0.5pt sym}}$. In \mycite{46} this polynomial
was denoted through $\tilde p_{\kern 1pt 1}$:
$$
\hskip -2em
\tilde p_{\kern 1pt 1}=p_{\kern 1pt 0}=x_1^2+x_2^2+x_3^2-L^2.
\mytag{1.11}
$$
In addition to \mythetag{1.11} in \mycite{46} the following seven polynomials
were considered:
$$
\gather
\hskip -2em
\aligned
\tilde p_{\kern 1pt 2}&=p_{\kern 1pt 1}+p_{\kern 1pt 2}+p_{\kern 1pt 3}
=(x_2^2+x_3^2-d_1^{\kern 1pt 2})\,+\\
&+\,(x_3^2+x_1^2-d_2^{\kern 1pt 2})+(x_1^2+x_2^2-d_3^{\kern 1pt 2}),
\endaligned\qquad\qquad
\mytag{1.12}\\
\vspace{2ex}
\aligned
\tilde p_{\kern 1pt 3}&=d_1\,p_{\kern 1pt 1}+d_2\,p_{\kern 1pt 2}+d_3
\,p_{\kern 1pt 3}=d_1\,(x_2^2+x_3^2-d_1^{\kern 1pt 2})\,+\\
&+\,d_2\,(x_3^2+x_1^2-d_2^{\kern 1pt 2})+d_3\,(x_1^2+x_2^2-d_3^{\kern 1pt 2}),
\endaligned\qquad\qquad
\mytag{1.13}\\
\vspace{2ex}
\aligned
\tilde p_{\kern 1pt 4}&=x_1\,p_{\kern 1pt 1}+x_2\,p_{\kern 1pt 2}+x_3
\,p_{\kern 1pt 3}=x_1\,(x_2^2+x_3^2-d_1^{\kern 1pt 2})\,+\\
&+\,x_2\,(x_3^2+x_1^2-d_2^{\kern 1pt 2})+x_3\,(x_1^2+x_2^2-d_3^{\kern 1pt 2}),
\endaligned\qquad\qquad
\mytag{1.14}\\
\vspace{2ex}
\aligned
\tilde p_{\kern 1pt 5}&=x_1\,d_1\,p_{\kern 1pt 1}+x_2\,d_2\,p_{\kern 1pt 2}
+x_3\,d_3\,p_{\kern 1pt 3}=x_1\,d_1\,(x_2^2+x_3^2-d_1^{\kern 1pt 2})\,+\\
&+\,x_2\,d_2\,(x_3^2+x_1^2-d_2^{\kern 1pt 2})+x_3\,d_3\,(x_1^2+x_2^2
-d_3^{\kern 1pt 2}),
\endaligned\qquad\qquad
\mytag{1.15}\\
\vspace{2ex}
\aligned
\tilde p_{\kern 1pt 6}&=x_1^2\,p_{\kern 1pt 1}+x_2^2\,p_{\kern 1pt 2}+x_3^2
\,p_{\kern 1pt 3}=x_1^2\,(x_2^2+x_3^2-d_1^{\kern 1pt 2})\,+\\
&+\,x_2^2\,(x_3^2+x_1^2-d_2^{\kern 1pt 2})+x_3^2\,(x_1^2+x_2^2
-d_3^{\kern 1pt 2}),
\endaligned\qquad\qquad
\mytag{1.16}\\
\vspace{2ex}
\aligned
\tilde p_{\kern 1pt 7}&=d_1^{\kern 1pt 2}\,p_{\kern 1pt 1}+d_2^{\kern 1pt 2}
\,p_{\kern 1pt 2}+d_3^{\kern 1pt 2}\,p_{\kern 1pt 3}=d_1^{\kern 1pt 2}
\,(x_2^2+x_3^2-d_1^{\kern 1pt 2})\,+\\
&+\,d_2^{\kern 1pt 2}\,(x_3^2+x_1^2-d_2^{\kern 1pt 2})+d_3^{\kern 1pt 2}
\,(x_1^2+x_2^2-d_3^{\kern 1pt 2}),
\endaligned\qquad\qquad
\mytag{1.17}\\
\vspace{2ex}
\aligned
\tilde p_{\kern 1pt 8}&=x_1^2\,d_1^{\kern 1pt 2}\,p_{\kern 1pt 1}+x_2^2
\,d_2^{\kern 1pt 2}\,p_{\kern 1pt 2}+x_3^2\,d_3^{\kern 1pt 2}
\,p_{\kern 1pt 3}=x_1^2\,d_1^{\kern 1pt 2}\,(x_2^2+x_3^2-d_1^{\kern 1pt 2})\,+\\
&+\,x_2^2\,d_2^{\kern 1pt 2}\,(x_3^2+x_1^2-d_2^{\kern 1pt 2})+x_3^2
\,d_3^{\kern 1pt 2}\,(x_1^2+x_2^2-d_3^{\kern 1pt 2}).
\endaligned\qquad\qquad
\mytag{1.18}
\endgather
$$
\mytheorem{1.1}The ideal $I_{\text{PC\kern -0.7pt\_\kern 0.5pt sym}}$ from
\mythetag{1.10} is finitely generated within the ring $\Sym\!\Bbb Q[M,L]$.
Eight polynomials \mythetag{1.11}, \mythetag{1.12}, \mythetag{1.13}, 
\mythetag{1.14}, \mythetag{1.15}, \mythetag{1.16}, \mythetag{1.17}, and 
\mythetag{1.18} belong to the ideal $I_{\text{PC\kern -0.7pt\_\kern 0.5pt sym}}$ 
and constitute a basis of this ideal.
\endproclaim
    The theorem~\mythetheorem{1.1} was proved in \mycite{46}. The ideal 
$I_{\text{C\kern -0.7pt\_\kern 0.5pt sym}}$ in \mythetag{1.8} is similar to the
ideal $I_{\text{PC\kern -0.7pt\_\kern 0.5pt sym}}$. There is the following 
theorem describing this ideal.  
\mytheorem{1.2}The ideal $I_{\text{C\kern -0.7pt\_\kern 0.5pt sym}}$ from
\mythetag{1.10} is finitely generated within the ring $\Sym\!\Bbb Q[M,L]$.
Seven polynomials \mythetag{1.12}, \mythetag{1.13}, \mythetag{1.14}, 
\mythetag{1.15}, \mythetag{1.16}, \mythetag{1.17}, and \mythetag{1.18} 
belong to the ideal $I_{\text{C\kern -0.7pt\_\kern 0.5pt sym}}$ 
and constitute a basis of this ideal.
\endproclaim
    The theorem~\mythetheorem{1.2} can be proved in a way similar to the
proof of the theorem~\mythetheorem{1.1} in \mycite{46}. I do not give the
proof of the theorem~\mythetheorem{1.2} here for the sake of brevity. 
\par
     Relying on the theorem~\mythetheorem{1.2} and using the polynomials
\mythetag{1.12}, \mythetag{1.13}, \mythetag{1.14}, \mythetag{1.15}, 
\mythetag{1.16}, \mythetag{1.17}, \mythetag{1.18}, we write the system 
of seven factor equations
$$
\xalignat 4
&\hskip -2em
&&\tilde p_{\kern 1pt 2}=0,&&\tilde p_{\kern 1pt 3}=0,&&\tilde p_{\kern 1pt 4}=0,
\qquad\\
\vspace{-1.7ex}
\mytag{1.19}\\
\vspace{-1.7ex}
&\hskip -2em
\tilde p_{\kern 1pt 5}=0,&&\tilde p_{\kern 1pt 6}=0,&&\tilde p_{\kern 1pt 7}=0,
&&\tilde p_{\kern 1pt 8}=0.\qquad
\endxalignat
$$
The factor equations \mythetag{1.19} correspond to the case of Euler cuboids.
Similarly, in the case of perfect cuboids, relying on the 
theorem~\mythetheorem{1.1} and using the polynomials given by the formulas  
\mythetag{1.11}, \mythetag{1.12}, \mythetag{1.13}, \mythetag{1.14}, 
\mythetag{1.15}, \mythetag{1.16}, \mythetag{1.17}, \mythetag{1.18}, we write 
the following system of eight factor equations:
$$
\xalignat 4
&\hskip -2em
\tilde p_{\kern 1pt 1}=0,
&&\tilde p_{\kern 1pt 2}=0,&&\tilde p_{\kern 1pt 3}=0,&&\tilde p_{\kern 1pt 4}=0,
\qquad\\
\vspace{-1.7ex}
\mytag{1.20}\\
\vspace{-1.7ex}
&\hskip -2em
\tilde p_{\kern 1pt 5}=0,&&\tilde p_{\kern 1pt 6}=0,&&\tilde p_{\kern 1pt 7}=0,
&&\tilde p_{\kern 1pt 8}=0.\qquad
\endxalignat
$$
The structure of the polynomials $\tilde p_{\kern 1pt 1}$, 
$\tilde p_{\kern 1pt 2}$, $\tilde p_{\kern 1pt 3}$, $\tilde p_{\kern 1pt 4}$, 
$\tilde p_{\kern 1pt 5}$, $\tilde p_{\kern 1pt 6}$, $\tilde p_{\kern 1pt 7}$, 
$\tilde p_{\kern 1pt 8}$ in \mythetag{1.11}, \mythetag{1.12}, \mythetag{1.13}, 
\mythetag{1.14}, \mythetag{1.15}, \mythetag{1.16}, \mythetag{1.17}, \mythetag{1.18}
is so that each solution of the equations \mythetag{1.2} is a solution for the 
equations \mythetag{1.19}. Similarly, each solution of the equations \mythetag{1.3} 
is a solution for the equations \mythetag{1.20}. The main goal of this paper is to 
prove converse propositions. They are given by the following two theorems.
\mytheorem{1.3} Each integer or rational solution of the factor equations 
\mythetag{1.19} such that $x_1>0$, $x_2>0$, $x_3>0$, $d_1>0$, $d_2>0$, 
and $d_3>0$ is an integer or rational solution for the equations 
\mythetag{1.2}.
\endproclaim
\mytheorem{1.4} Each integer or rational solution of the factor equations 
\mythetag{1.20} such that $x_1>0$, $x_2>0$, $x_3>0$, $d_1>0$, $d_2>0$, 
and $d_3>0$ is an integer or rational solution for the equations 
\mythetag{1.3}.
\endproclaim
\head
2. The analysis of the factor equations.
\endhead
      Let's consider the factor equations \mythetag{1.19} associated with Euler
cuboids. Due to \mythetag{1.12}, \mythetag{1.13}, \mythetag{1.14}, 
\mythetag{1.15}, \mythetag{1.16}, \mythetag{1.17}, and \mythetag{1.18} the
factor equations \mythetag{1.19} can be united into a single matrix equation
$$
\hskip -2em
\Vmatrix
1 & 1 & 1\\
d_1 & d_2 & d_3\\
x_1 & x_2 & x_3\\
x_1\,d_1 & x_2\,d_2 & x_3\,d_3\\
\vspace{0.5ex}
x_1^2 & x_2^2 & x_3^2\\
\vspace{0.7ex}
d_1^{\kern 1pt 2} & d_2^{\kern 1pt 2} & d_3^{\kern 1pt 2}\\
\vspace{0.7ex}
x_1^2\,d_1^{\kern 1pt 2} & x_2^2\,d_2^{\kern 1pt 2} & x_3^2\,d_3^{\kern 1pt 2}
\endVmatrix
\cdot\Vmatrix
p_{\kern 1pt 1}\\
p_{\kern 1pt 2}\\
p_{\kern 1pt 3}
\endVmatrix=\Vmatrix
0\\ 0\\ 0\\ 0\\ 0\\ 0\\ 0
\endVmatrix.
\mytag{2.1}
$$
In order to study the equations \mythetag{2.1} we denote through $N$ the 
transposed matrix 
$$
\hskip -2em
N=\Vmatrix
1 & d_1 & x_1 & x_1\,d_1 & x_1^2 & d_1^{\kern 1pt 2} & x_1^2\,d_1^{\kern 1pt 2}\\
\vspace{1ex}
1 & d_2 & x_2 & x_2\,d_2 & x_2^2 & d_2^{\kern 1pt 2} & x_2^2\,d_2^{\kern 1pt 2}\\
\vspace{1ex}
1 & d_3 & x_3 & x_3\,d_3 & x_3^2 & d_3^{\kern 1pt 2} & x_3^2\,d_3^{\kern 1pt 2}
\endVmatrix.
\mytag{2.2}
$$
If we have a solution of the equation \mythetag{2.1} which is not a solution 
for the initial system of cuboid equations \mythetag{1.2}, then the equations 
\mythetag{1.2} should not be fulfilled simultaneously. Therefore we have the 
vectorial inequality 
$$
\hskip -2em
\Vmatrix
p_{\kern 1pt 1}\\
p_{\kern 1pt 2}\\
p_{\kern 1pt 3}
\endVmatrix\neq 0.
\mytag{2.3}
$$
Applying \mythetag{2.3} to \mythetag{2.1}, we derive that the columns 
of the matrix in \mythetag{2.1} are linearly dependent. Then the rows of $N$ in 
\mythetag{2.2} are also linearly dependent, i\.\,e\.
$$
\hskip -2em
\rank N\leqslant 2.
\mytag{2.4}
$$
The condition \mythetag{2.4} leads to several special cases which are 
considered below one by one. In addition to $N$ we define the following
two matrices:
$$
\xalignat 2
&\hskip -2em
N_1=\Vmatrix
1 & d_1\\
\vspace{0.5ex}
1 & d_2\\
\vspace{0.5ex}
1 & d_3
\endVmatrix,
&&N_2=\Vmatrix
1 & x_1\\
\vspace{0.5ex}
1 & x_2\\
\vspace{0.5ex}
1 & x_3
\endVmatrix.
\mytag{2.5}
\endxalignat
$$
The matrices $N_1$ and $N_2$ in \mythetag{2.5} are submatrices of the matrix $N$.
\par
\head
3. The case $\rank N=1$.
\endhead
     The first column of the matrix \mythetag{2.2} is nonzero. Therefore 
$\rank N>0$. Now we consider the case where $\rank N=1$. In this case 
each column of the matrix $N$ is proportional to its first column. 
In particular, this yields
$$
\xalignat 2
&\hskip -2em
\Vmatrix x_1\\ x_2\\ x_3\endVmatrix
=\alpha\cdot\Vmatrix 1\\ 1\\ 1\endVmatrix,
&&\Vmatrix d_1\\ d_2\\ d_3\endVmatrix
=\beta\cdot\Vmatrix 1\\ 1\\ 1\endVmatrix.
\mytag{3.1}
\endxalignat
$$
The equations \mythetag{3.1} lead to the equalities 
$$
\xalignat 2
&\hskip -2em
x_1=x_2=x_3,
&&d_1=d_2=d_3.
\mytag{3.2}
\endxalignat
$$
Applying \mythetag{3.2} to the formulas \mythetag{1.1}, we derive 
$$
\hskip -2em
p_{\kern 1pt 1}=p_{\kern 1pt 2}=p_{\kern 1pt 3}.
\mytag{3.3}
$$
Then we substitute \mythetag{3.3} into \mythetag{1.12}. As a result we get 
$$
\hskip -2em
\tilde p_{\kern 1pt 2}=3\,p_{\kern 1pt 1}=3\,p_{\kern 1pt 2}
=3\,p_{\kern 1pt 3}.
\mytag{3.4}
$$
The relationships \mythetag{3.4} mean that if the equations \mythetag{1.19}
are fulfilled, then in the case of $\rank N=1$ the equations \mythetag{1.2}
are also fulfilled. 
\mytheorem{3.1} Each solution of the equations \mythetag{1.19} corresponding 
to the case $\rank N=1$ is a solution for the equations \mythetag{1.2}. 
\endproclaim
\mytheorem{3.2} Each solution of the equations \mythetag{1.20} corresponding 
to the case $\rank N=1$ is a solution for the equations \mythetag{1.3}. 
\endproclaim
Due to \mythetag{1.11} the equation $p_{\kern 1pt 0}=0$ in \mythetag{1.3} 
coincides with the equation $\tilde p_{\kern 1pt 1}=0$ in \mythetag{1.20}.
For this reason the theorem~\mythetheorem{3.2} is immediate from the
theorem~\mythetheorem{3.1}.
\head
4. The case $\rank N_1=2$ and $\rank N_2=1$.
\endhead
    The condition $\rank N_2=1$ for the matrix $N_2$ in \mythetag{2.5} means
that the third column of the matrix \mythetag{2.2} is proportional to the
first column of this matrix. \pagebreak The condition $\rank N_1=2$ for the matrix $N_1$ 
in \mythetag{2.5} means that the first and the second columns of the matrix 
\mythetag{2.2} are linearly independent. Other columns are expressed as linear
combinations of these two columns. As a result we can write
$$
\xalignat 2
&\hskip -2em
\Vmatrix x_1\\ x_2\\ x_3\endVmatrix
=\alpha\cdot\Vmatrix 1\\ 1\\ 1\endVmatrix,
&&\Vmatrix d_1^{\kern 1pt 2}\\ d_2^{\kern 1pt 2}\\ d_3^{\kern 1pt 2}
\endVmatrix=\beta\cdot\Vmatrix d_1\\ d_2\\ d_3\endVmatrix
+\gamma\cdot\Vmatrix 1\\ 1\\ 1\endVmatrix.
\mytag{4.1}
\endxalignat
$$
It is easy to see that the conditions \mythetag{4.1} are sufficient for
to provide the condition $\rank N=2$, which is in agreement with 
\mythetag{2.4}.\par
     The second equality in \mythetag{4.1} is very important. It means
that $d_1$, $d_2$, and $d_3$, are roots of the following quadratic equation:
$$
\hskip -2em
d^{\kern 1pt 2}-\beta\,d-\gamma=0.
\mytag{4.2}
$$
The quadratic equation \mythetag{4.2} has at most two roots. Let's denote 
them $s_1$ and $s_2$. Then we have the following subcases derived from 
$\rank N_1=2$ and $\rank N_2=1$:
$$
\xalignat 3
&\hskip -2em
\Vmatrix d_1\\ d_2\\ d_3\endVmatrix
=\Vmatrix s_1\\ s_1\\ s_2\endVmatrix,
&&\Vmatrix d_1\\ d_2\\ d_3\endVmatrix
=\Vmatrix s_1\\ s_2\\ s_1\endVmatrix,
&&\Vmatrix d_1\\ d_2\\ d_3\endVmatrix
=\Vmatrix s_2\\ s_1\\ s_1\endVmatrix.
\quad
\mytag{4.3}
\endxalignat
$$
The numbers $s_1$ and $s_2$ in the formulas \mythetag{4.3} are arbitrary 
two numbers not coinciding with each other: $s_1\neq s_2$. They are integer
numbers in the case of integer solutions and they are rational numbers in
the case of rational solutions.\par
     The three cases in \mythetag{4.3} are similar to each other. Without loss
of generality we can consider only one of them, e\.\,g\. the first one. Then 
from \mythetag{4.3} we derive 
$$
\gather
\hskip -2em
(s_1-s_2)\cdot
\Vmatrix 1\\ 1\\ 0\endVmatrix
=\Vmatrix s_1\\ s_1\\ s_2\endVmatrix
-s_2\cdot\Vmatrix 1\\ 1\\ 1\endVmatrix,
\mytag{4.4}\\
\vspace{2ex}
\hskip -2em
(s_2-s_1)\cdot
\Vmatrix 0\\ 0\\ 1\endVmatrix
=\Vmatrix s_1\\ s_1\\ s_2\endVmatrix
-s_1\cdot\Vmatrix 1\\ 1\\ 1\endVmatrix.
\mytag{4.5}
\endgather
$$
Due to the relationships \mythetag{4.4} and \mythetag{4.5} the matrix equation 
\mythetag{2.1} reduces to 
$$
\hskip -2em
\Vmatrix
1 & 1 & 0\\
0 & 0 & 1
\endVmatrix
\cdot\Vmatrix
p_{\kern 1pt 1}\\
p_{\kern 1pt 2}\\
p_{\kern 1pt 3}
\endVmatrix=\Vmatrix
0\\ 0
\endVmatrix.
\mytag{4.6}
$$
The matrix equality \mythetag{4.6} means that instead of the seven equations 
\mythetag{1.19} we have two equations $p_{\kern 1pt 1}+p_{\kern 1pt 2}=0$ 
and $p_{\kern 1pt 3}=0$. Substituting $x_1=x_2=x_3=\alpha$, $d_1=d_2=s_1$,
and $d_3=s_2$ into these two equations, we derive 
$$
\xalignat 2
&\hskip -2em
s_1^2-2\,\alpha^2=0,
&&s_2^2-2\,\alpha^2=0.
\mytag{4.7}
\endxalignat
$$
The equations \mythetag{4.7} can be written in the following way:
$$
\xalignat 2
&\hskip -2em
|s_1|=\sqrt{2}\,|\alpha|,
&&|s_2|=\sqrt{2}\,|\alpha|.
\mytag{4.8}
\endxalignat
$$ 
Now it is easy to see that the equations \mythetag{4.8} can be satisfied by
three integer or rational numbers $s_1$, $s_2$, and $\alpha$ if and only if all 
of them are zero. Substituting $s_1=s_2=\alpha=0$ into \mythetag{4.1} and
\mythetag{4.3}, we get 
$$
\xalignat 2
&\hskip -2em
x_1=x_2=x_3=0,
&&d_1=d_2=d_3=0.
\mytag{4.9}
\endxalignat
$$ 
The equalities \mythetag{4.9} contradict the condition $\rank N_1=2$ for
the matrix $N_1$ in \mythetag{2.5}. This contradiction yields the following
two theorems. 
\mytheorem{4.1} The factor equations \mythetag{1.19}, as well as the original
equations \mythetag{1.2}, have no integer or rational solution in the case of 
$\rank N_1=2$ and $\rank N_2=1$.
\endproclaim
\mytheorem{4.2} The factor equations \mythetag{1.20}, as well as the original
equations \mythetag{1.3}, have no integer or rational solution in the case of 
$\rank N_1=2$ and $\rank N_2=1$.
\endproclaim
\head
5. The case $\rank N_1=1$ and $\rank N_2=2$.
\endhead
    The condition $\rank N_1=1$ for the matrix $N_1$ in \mythetag{2.5} means
that the second column of the matrix \mythetag{2.2} is proportional to the
first column of this matrix. The condition $\rank N_2=2$ for the matrix $N_2$ 
in \mythetag{2.5} means that the first and the third columns of the matrix 
\mythetag{2.2} are linearly independent. Other columns are expressed as linear
combinations of these two columns. As a result we can write the 
relationships similar to the relationships \mythetag{4.1}: 
$$
\xalignat 2
&\hskip -2em
\Vmatrix d_1\\ d_2\\ d_3\endVmatrix
=\delta\cdot\Vmatrix 1\\ 1\\ 1\endVmatrix,
&&\Vmatrix x_1^2\\ x_2^2\\ x_3^2
\endVmatrix=\varepsilon\cdot\Vmatrix x_1\\ x_2\\ x_3\endVmatrix
+\zeta\cdot\Vmatrix 1\\ 1\\ 1\endVmatrix.
\mytag{5.1}
\endxalignat
$$
The conditions \mythetag{5.1} are sufficient for to provide the condition 
$\rank N=2$.\par
     Like in the case of \mythetag{4.1}, the second condition \mythetag{5.1}
mean that $x_1$, $x_2$, and $x_3$ are roots of the quadratic equation
similar to \mythetag{4.2}:
$$
\hskip -2em
x^2-\varepsilon\,x-\zeta=0.
\mytag{5.2}
$$
The quadratic equation \mythetag{5.2} has at most two roots. Let's denote 
them $r_1$ and $r_2$. Then we have the following subcases derived from 
$\rank N_1=1$ and $\rank N_2=2$:
$$
\xalignat 3
&\hskip -2em
\Vmatrix x_1\\ x_2\\ x_3\endVmatrix
=\Vmatrix r_1\\ r_1\\ r_2\endVmatrix,
&&\Vmatrix x_1\\ x_2\\ x_3\endVmatrix
=\Vmatrix r_1\\ r_2\\ r_1\endVmatrix,
&&\Vmatrix x_1\\ x_2\\ x_3\endVmatrix
=\Vmatrix r_2\\ r_1\\ r_1\endVmatrix.
\quad
\mytag{5.3}
\endxalignat
$$
The numbers $r_1$ and $r_2$ in the formulas \mythetag{5.3} are arbitrary two 
integer or rational numbers not coinciding with each other: $r_1\neq r_2$.
\par
     The three cases in \mythetag{5.3} are similar to each other. Without loss
of generality we can consider only one of them, e\.\,g\. the first one. Then 
from \mythetag{5.3} we derive      
$$
\allowdisplaybreaks
\gather
\hskip -2em
(r_1-r_2)\cdot\Vmatrix 1\\ 1\\ 0\endVmatrix
=\Vmatrix r_1\\ r_1\\ r_2\endVmatrix
-r_2\cdot\Vmatrix 1\\ 1\\ 1\endVmatrix,
\mytag{5.4}\\
\vspace{2ex}
\hskip -2em
(r_2-r_1)\cdot\Vmatrix 0\\ 0\\ 1\endVmatrix
=\Vmatrix r_1\\ r_1\\ r_2\endVmatrix
-r_1\cdot\Vmatrix 1\\ 1\\ 1\endVmatrix.
\mytag{5.5}
\endgather
$$
Due to the relationships \mythetag{5.4} and \mythetag{5.5} the matrix equation 
\mythetag{2.1} reduces to 
$$
\hskip -2em
\Vmatrix
1 & 1 & 0\\
0 & 0 & 1
\endVmatrix
\cdot\Vmatrix
p_{\kern 1pt 1}\\
p_{\kern 1pt 2}\\
p_{\kern 1pt 3}
\endVmatrix=\Vmatrix
0\\ 0
\endVmatrix.
\mytag{5.6}
$$
The matrix equality \mythetag{5.6} means that instead of the seven equations 
\mythetag{1.19} we have two equations $p_{\kern 1pt 1}+p_{\kern 1pt 2}=0$ 
and $p_{\kern 1pt 3}=0$. Substituting $d_1=d_2=d_3=\delta$, $x_1=x_2=r_1$,
and $x_3=r_2$ into these two equations, we derive 
$$
\xalignat 2
&\hskip -2em
r_1^2+r_2^2-\delta^2=0,
&&2\,r_1^2-\delta^2=0.
\quad
\mytag{5.7}
\endxalignat
$$
The second equation \mythetag{5.7} can be written in the following form:
$$
\hskip -2em
|\delta|=\sqrt{2}\,|r_1|.
\mytag{5.8}
$$ 
The equation \mythetag{5.8} can be satisfied by two integer or rational 
numbers $r_1$ and $\delta$ if and only if both of them are zero. 
Substituting $r_1=\delta=0$ into the first equation \mythetag{5.7}, we get 
$r_2=0$. Substituting $r_1=r_2=\delta=0$ into \mythetag{5.1} and
\mythetag{5.3}, we get 
$$
\xalignat 2
&\hskip -2em
x_1=x_2=x_3=0,
&&d_1=d_2=d_3=0.
\mytag{5.9}
\endxalignat
$$ 
The equalities \mythetag{5.9} contradict the condition $\rank N_2=2$ for
the matrix $N_2$ in \mythetag{2.5}. This contradiction yields the following
two theorems.
\mytheorem{5.1} The factor equations \mythetag{1.19}, as well as the original
equations \mythetag{1.2}, have no integer or rational solution in the case of 
$\rank N_1=1$ and $\rank N_2=2$.
\endproclaim
\mytheorem{5.2} The factor equations \mythetag{1.20}, as well as the original
equations \mythetag{1.3}, have no integer or rational solution in the case of 
$\rank N_1=1$ and $\rank N_2=2$.
\endproclaim
\head
6. The case $\rank N_1=2$ and $\rank N_2=2$.
\endhead
    In this case the columns of both matrices $N_1$ and $N_2$ in \mythetag{2.5}
are linearly independent. Hence each column of the matrix \mythetag{2.2} can be
expressed as a linear combination of the first and the second columns of this 
matrix or as a linear combination of the first and the third columns of this 
matrix. In particular, we have
$$
\xalignat 2
&\hskip -2em
\Vmatrix d_1^{\kern 1pt 2}\\ d_2^{\kern 1pt 2}\\ d_3^{\kern 1pt 2}
\endVmatrix=\beta\cdot\Vmatrix d_1\\ d_2\\ d_3\endVmatrix
+\gamma\cdot\Vmatrix 1\\ 1\\ 1\endVmatrix,
&&\Vmatrix x_1^2\\ x_2^2\\ x_3^2
\endVmatrix=\varepsilon\cdot\Vmatrix x_1\\ x_2\\ x_3\endVmatrix
+\zeta\cdot\Vmatrix 1\\ 1\\ 1\endVmatrix.
\mytag{6.1}
\endxalignat
$$
The relationships \mythetag{6.1} mean that $x_1$, $x_2$, $x_3$ and $d_1$, $d_2$, 
$d_3$ are roots of two quadratic equations coinciding with \mythetag{5.2} and 
\mythetag{4.2} respectively. As a result we distinguish three subcases 
\mythetag{4.3} with $s_1\neq s_2$ and three subcases \mythetag{5.3} with $r_1
\neq r_2$. The first subcase \mythetag{4.3} should be paired with the first 
subcase \mythetag{5.3}, the second subcase \mythetag{4.3} should be paired with 
the second subcase \mythetag{5.3}, and the third subcase \mythetag{4.3} should 
be paired with the third subcase \mythetag{5.3}. Otherwise we would have $\rank N
\geqslant 3$, which contradicts the condition \mythetag{2.4}.\par
     Due to the pairing of subcases we have three subcases instead of nine ones,
which are a priori possible. These three subcases are similar to each other. 
Therefore without loss of generality we can consider only one subcase, e\.\,g\. 
the following one: 
$$
\xalignat 2
&\hskip -2em
\Vmatrix d_1\\ d_2\\ d_3\endVmatrix
=\Vmatrix s_1\\ s_1\\ s_2\endVmatrix,
&&\Vmatrix x_1\\ x_2\\ x_3\endVmatrix
=\Vmatrix r_1\\ r_1\\ r_2\endVmatrix.
\mytag{6.2}
\endxalignat
$$
Here $s_1\neq s_2$ and $r_1\neq r_2$. The relationships \mythetag{6.2} lead to
the relationships \mythetag{4.4}, \mythetag{4.5}, \mythetag{5.4}, \mythetag{5.5} 
and then to the equations \mythetag{4.6} and \mythetag{5.6}. The matrix equations 
\mythetag{4.6} and \mythetag{5.6} mean that instead of the seven equations 
\mythetag{1.19} we have two equations $p_{\kern 1pt 1}+p_{\kern 1pt 2}=0$ and 
$p_{\kern 1pt 3}=0$. Substituting $d_1=d_2=s_1$, $d_3=s_2$, $x_1=x_2=r_1$, and 
$x_3=r_2$ into these two equations, we derive 
$$
\xalignat 2
&\hskip -2em
r_1^2+r_2^2-s_1^2=0,
&&2\,r_1^2-s_2^2=0. 
\mytag{6.3}
\endxalignat
$$
The second equation \mythetag{6.3} can be written in the following way:
$$
\hskip -2em
\sqrt{2}\,|r_1|=|s_2|. 
\mytag{6.4}
$$
The equation \mythetag{6.4} can be satisfied by two integer or rational 
numbers $r_1$ and $s_2$ if and only if both of them are zero. Substituting
$r_1=s_2=0$ into \mythetag{6.3}, we get $|r_2|=|s_1|=\theta$. Then the
equalities \mythetag{6.2} are written as 
$$
\xalignat 2
&\hskip -2em
\Vmatrix d_1\\ d_2\\ d_3\endVmatrix
=\pm\Vmatrix \theta\\ \theta\\ 0\endVmatrix,
&&\Vmatrix x_1\\ x_2\\ x_3\endVmatrix
=\pm\Vmatrix 0\\ 0\\ \theta\endVmatrix.
\mytag{6.5}
\endxalignat
$$
The equalities \mythetag{6.5} lead to the equalities $x_1=x_2=0$ and $d_3=0$. 
The latter ones contradict the inequalities in the theorems~\mythetheorem{1.3}
and \mythetheorem{1.4}. Therefore we can conclude this section with the 
following two theorems. 
\mytheorem{6.1} The factor equations \mythetag{1.19}, as well as the original 
equations \mythetag{1.2}, have no integer or rational solutions such that  
$x_1>0$, $x_2>0$, $x_3>0$, $d_1>0$, $d_2>0$, and $d_3>0$ in the case of 
$\rank N_1=2$ and $\rank N_2=2$.
\endproclaim
\mytheorem{6.2} The factor equations \mythetag{1.20}, as well as the original 
equations \mythetag{1.3}, have no integer or rational solutions such that  
$x_1>0$, $x_2>0$, $x_3>0$, $d_1>0$, $d_2>0$, and $d_3>0$ in the case of 
$\rank N_1=2$ and $\rank N_2=2$.
\endproclaim
\head
7. The ultimate result and conclusions. 
\endhead
     The four cases considered in sections 3, 4, 5, and 6 exhaust all 
options compatible with the inequality \mythetag{2.4}. For this reason the 
theorems~\mythetheorem{1.3} and \mythetheorem{1.4} follow from the 
theorems~\mythetheorem{3.1}, \mythetheorem{4.1}, \mythetheorem{5.1}, 
\mythetheorem{6.1} and the theorems~\mythetheorem{3.2}, \mythetheorem{4.2}, 
\mythetheorem{5.2}, \mythetheorem{6.2} respectively. The 
theorems~\mythetheorem{1.3} and \mythetheorem{1.4} constitute the main result 
of this paper. \pagebreak The theorem~\mythetheorem{1.4} means that the 
factor equations \mythetag{1.20} are equally admissible for seeking perfect 
cuboids or for proving their non-existence as the original equations 
\mythetag{1.3}. As for the factor equations \mythetag{1.19}, due to the 
theorem~\mythetheorem{1.3} they are equally admissible for selecting Euler 
cuboids as the original equations \mythetag{1.2}.

\enddocument
\end